# ABOUT VERY PERFECT NUMBERS

Mihály Bencze and Florin Popovici
Department of Mathematics
Áprily Lajos College
Braşov, Romania

Florentin Smarandache, Ph D
Professor
Chair of Department of Math & Sciences
University of New Mexico
200 College Road
Gallup, NM 87301, USA
E-mail: smarand@unm.edu

**Abstract.**
In this short paper we prove that the square of an odd prime number cannot be a very perfect number.

**Introduction.**
A natural number $n$ is called *very perfect* if $\sigma(\sigma(n)) = 2n$ (see [1]), where $\sigma(x)$ means the sum of all positive divisors of the natural number x.

We now prove the following result:
**Theorem.** The square of an odd prime number cannot be a very perfect number.

*Proof:* Let's consider $n = p^2$, where $p$ is an odd prime number, then
$$\sigma(n) = 1 + p + p^2, \quad \sigma(\sigma(n)) = \sigma(1 + p + p^2) = 2p^2.$$
We decompose $\sigma(n)$ in canonical form, from where $1 + p + p^2 = p_1^{\alpha_1} p_2^{\alpha_2} ... p_k^{\alpha_k}$. Because $p(p+1)+1$ is odd, in the canonical decompose there must be only odd numbers.
$$\sigma(\sigma(n)) = \left(1 + p_1 + ... + p_1^{\alpha}\right)...\left(1 + p_k + ... + p_k^{\alpha_k}\right) = \frac{p_1^{\alpha_1+1} - 1}{p_1 - 1} ... \frac{p_k^{\alpha_k+1} - 1}{p_k - 1} = 2p^2$$
Because
$$\frac{p_1^{\alpha_1+1} - 1}{p_1 - 1} > 2, ..., \frac{p_k^{\alpha_k+1} - 1}{p_k - 1} > 2$$
one obtains that $2p^2$ cannot be decomposed in more than two factors, then each one is $> 2$, therefore $k \le 2$.
**Case 1.** For $k = 1$ we find $\sigma(n) = 1 + p + p^2 = p_1^{\alpha_1}$, from where one obtains

$$p_1^{\alpha_1+1} = p_1(1+p+p^2) \text{ and}$$

$$\sigma(\sigma(n)) = \frac{p_1^{\alpha_1+1}-1}{p_1-1} = 2p^2, \ p_1(1+p+p^2)-1 = 2p^2(p_1-1),$$

from where

$$p_1 - 1 = p(pp_1 - 2p - p_1).$$

The right side is divisible by $p$, thus $p_1 - 1$ is a $p$ multiple. Because $p_1 > 2$ it results that

$$p_1 \geq p-1 \text{ and } p_1^2 \geq (p+1)^2 > p^2+p+1 = p_1^{\alpha_1},$$

thus $\alpha_1 = 1$ and

$$\sigma(n) = p^2 + p + 1 = p_1, \ \sigma(\sigma(n)) = \sigma(p_1) = 1 + p_1.$$

If $n$ is very perfect then $1 + p_1 = 2p^2$ or $p^2 + p + 2 = 2p^2$. The solutions of the equation are $p = -1$, and $p = 2$ which is a contradiction.

**Case 2.** For $k = 2$ we have $\sigma(n) = p^2 + p + 1 = p_1^{\alpha_1} p_2^{\alpha_2}$.

$$\sigma(\sigma(n)) = (1+p_1+...+p_1^{\alpha_1})(1+p_2+...+p_2^{\alpha_2}) = \frac{p_1^{\alpha_1+1}-1}{p_1-1} \cdot \frac{p_2^{\alpha_2+1}-1}{p_2-1} = 2p^2.$$

Because

$$\frac{p_1^{\alpha_1+1}-1}{p_1-1} > 2 \text{ and } \frac{p_2^{\alpha_2+1}-1}{p_2-1} > 2,$$

it results

$$\frac{p_1^{\alpha_1+1}}{p_1-1} = p \text{ and } \frac{p_2^{\alpha_2+1}-1}{p_2-1} = 2p$$

(or inverse), thus

$$p_1^{\alpha_1+1} - 1 = p(p_1-1), \ p_2^{\alpha_2+1} - 1 = 2p(p_2-1),$$

then

$$p_1^{\alpha_1+1} p_2^{\alpha_2+1} - p_1^{\alpha_1+1} - p_2^{\alpha_2+1} + 1 = 2p^2(p_1-1)(p_2-1),$$

thus

$$\sigma(n) = p^2 + p + 1 = p_1^{\alpha_1+1} p_2^{\alpha_2+1}$$

and

$$p_1 p_2 (p^2+p+1) = 2p^2(p_1-1)(p_2-1) + p_1^{\alpha_1+1} + p_2^{\alpha_2+1} - 1$$

or

$$p_1 p_2 p(p+1) + p_1 p_2 - 1 = 2p^2(p_1-1)(p_2-1) + (p_1^{\alpha_1+1}-1) + (p_2^{\alpha_2+1}-1) =$$
$$= 2p^2(p_1-1)(p_2-1) + p(p_1-1) + 2p(p_2-1)$$

accordingly $p$ divides $p_1 p_2 - 1$, thus $p_1 p_2 > p+1$ and

$$p_1^2 p_2^2 \geq (p+1)^2 > p^2 + p + 1 = p_1^{\alpha_1} p_2^{\alpha_2}.$$

Hence:

$\Pi_1$) If $\alpha_1 = 1$ and $n = 2p^2$,

then
$$\sigma(n) = p^2 + p + 1 = p_1 p_2^{\alpha_2} \text{ and } \frac{p_1^2 - 1}{p_1 - 1} = p, \text{ and } \frac{p_2^{\alpha_2+1} - 1}{p_2 - 1} = 2p,$$
thus $p_1 + 1 = p$ which is a contradiction.

$\Pi_2$ ) If $\alpha_2 = 1$ and $n = 2p^2$,
then
$$\sigma(n) = p^2 + p + 1 = p_1^{\alpha_1} p_2 \text{ and } \frac{p_1^{\alpha_1+1} - 1}{p_1 - 1} = p, \text{ and}$$

$$\frac{p_2^2 - 1}{p_2 - 1} = 2p,$$

thus
$$p_2 + 1 = 2p, \ p_2 = 2p - 1$$
and
$$\sigma(n) = p^2 + p + 1 = p_1^{\alpha_1}(2p + 1),$$
from where
$$4\sigma(n) = (2p-1)(2p+3) + 7 = 4 p_1^{\alpha_1}(2p-1),$$
accordingly 7 is divisible by $2p - 1$ and thus $p$ is divisible by 4 which is a contradiction.

**Reference:**

[1] Suryanarayama – Elemente der Mathematik – 1969.